%% file: coupled.tex
\documentclass[a4paper]{coupled2021}

\usepackage{graphicx}
\usepackage{amsmath}
\usepackage{amsfonts}
\usepackage{amssymb}
\usepackage{mathtools}
\usepackage{xcolor}
\usepackage{graphicx}
\graphicspath{{Illustrations/}}
\usepackage{multirow}
\usepackage{booktabs}
\usepackage{csquotes}
\usepackage[english]{babel}
\usepackage{pifont} 
\usepackage{url}

\usepackage[
backend=biber,
bibencoding=utf8,
citestyle=numeric-comp,
giveninits=true,
doi=false,
isbn=false,
url=false,
maxcitenames=6,
maxbibnames=100,
]{biblatex}
\AtEveryBibitem{
  \clearfield{note}%
}
\DeclareNameAlias{author}{last-first} 
\DeclareFieldFormat[article, inbook, incollection, inproceedings, misc, thesis, unpublished]
{pages}{#1}

\renewbibmacro{in:}{}
\DeclareFieldFormat[inbook, book, incollection, inproceedings, misc, thesis, unpublished, report]
{date}{\mkbibparens{#1}}

\DeclareFieldFormat[article]
{volume}{\textbf{#1}}
\DeclareFieldFormat[article, inbook, incollection, inproceedings, misc, thesis, unpublished, report]
{number}{\mkbibparens{#1}}

\renewbibmacro*{volume+number+eid}{
  \printfield{volume}%
  \printfield{number}%
  \printunit{\addcolon}%
}

\renewbibmacro*{journal+issuetitle}{%
  \usebibmacro{journal}%
  \iffieldundef{series}%
    {}%
    {\printfield{series}%
     \setunit{\addspace}}%
  \usebibmacro{issue+date}%
  \setunit{\addspace}%
  \usebibmacro{issue}%
  \usebibmacro{volume+number+eid}%
  \newunit}

\newcommand{\ilu}[1]{\text{ILU(#1)}}
\newcommand{\dumux}{$\text{\textsc{DuMu}}^\text{\textsc{x}}$\ }
\newcommand{\dumuX}{$\text{\textsc{DuMu}}^\text{\textsc{x}}$}

\newcommand{\dunE}{\textsc{Dune}}

\newcommand{\istL}{\textsc{Istl}}
\newcommand{\umfpack}{\textsc{Umfpack}\ }
\newcommand{\umfpacK}{\textsc{Umfpack}}

\newcommand{\solv}{\mathcal{S}}
\newcommand{\postproc}{\mathcal{I}^\mathrm{post}}

\newcommand{\ff}{^\mathrm{ff}}
\newcommand{\prm}{^\mathrm{pm}}
\newcommand{\tp}{^\mathrm{td}}
\newcommand{\pv}{^\mathrm{pv}}
\newcommand{\pc}{\mathcal{P}}
\newcommand{\bgs}{_\mathrm{BGS}}
\newcommand{\bjac}{_\mathrm{BJ}}
\newcommand{\app}{_{A'}}
\newcommand{\dpp}{_{D'}}
\newcommand{\vu}{_{V}}

\newcommand{\yes}{\text{\ding{51}}}%
\newcommand{\no}{\text{\ding{55}}}%

\newcommand{\executeiffilenewer}[3]{\ifnum\pdfstrcmp{\pdffilemoddate{#1}}%
      {\pdffilemoddate{#2}}>0%
      {\immediate\write18{#3}}\fi%
}

\newcommand\blfootnote[1]{%
  \begingroup
  \renewcommand\thefootnote{}\footnote{#1}%
  \addtocounter{footnote}{-1}%
  \endgroup
}

\title{PARTITIONED COUPLING VS. MONOLITHIC BLOCK-PRECONDITIONING APPROACHES FOR SOLVING  STOKES-DARCY SYSTEMS}

\author{JENNY SCHMALFUSS$^{* \star}$, CEDRIC RIETHM{\"U}LLER$^{\dag \star}$, MIRCO ALTENBERND$^{\dag}$, KILIAN WEISHAUPT$^{\dag\dag}$ AND DOMINIK G{\"O}DDEKE$^{\dag\#}$}

\heading{Jenny Schmalfuss, Cedric Riethm{\"u}ller, Mirco Altenbernd, Kilian Weishaupt and Dominik G{\"o}ddeke}

\address{$^{*}$Institute for Visualization and Interactive Systems (VIS), University of Stuttgart\\
70569 Stuttgart, Germany. E-mail: jenny.schmalfuss@vis.uni-stuttgart.de
\and
$^{\dag}$Institute of Applied Analysis and Numerical Simulation (IANS), University of Stuttgart
\and
$^{\dag\dag}$Institute for Modelling Hydraulic and Environmental Systems (IWS), University of Stuttgart
\and
$^{\#}$Stuttgart Center for Simulation Science (SC SimTech), University of Stuttgart
}


\abstract{We consider the time-dependent Stokes-Darcy problem as a model case for the challenges involved in solving coupled systems. Keeping the model, its discretization, and the underlying numerics for the subproblems in the free-flow domain and the porous medium constant, we focus on different solver approaches for the coupled problem. We compare a partitioned coupling approach using the coupling library preCICE with a monolithic block-preconditioned one that is tailored to different formulations of the problem. Both approaches enable the reuse of already available iterative solvers and preconditioners, in our case, from the \dumuX{} framework.
Our results indicate that the approaches can yield performance and scalability improvements compared to using direct solvers: Partitioned coupling is able to solve large problems faster if iterative solvers with suitable preconditioners are applied for the subproblems. The monolithic approach shows even stronger requirements on preconditioning, as standard simple solvers fail to converge. Our monolithic block preconditioning yields the fastest runtimes for large systems, but they vary strongly with the preconditioner configuration. Interestingly, using a specialized Uzawa preconditioner for the Stokes subsystem leads to overall increased runtimes compared to block preconditioners utilizing a more general algebraic multigrid. This highlights that optimizing for the non-coupled cases does not always pay off.}

\begin{document}
\blfootnote{$\star$ The first two authors contributed equally to this paper.}

%
%
\section{INTRODUCTION}

Coupled systems of free flow adjacent to permeable media have a decisive role in many applications.
Examples include the environmental sciences (soil water evaporation),  medical contexts (intervascular exchange), material design (optimization of fuel cell water management) or technical applications (drying of perishable goods) to name just a few.
Capturing the complex interplay between the two flow domains is essential, however, the governing systems of equations form a coupled problem which can become quite complex to solve.
This even holds for single-phase-flow systems, such as a river flowing over its porous bed.
In this paper, we deliberately restrict ourselves to a simple, stationary, single-phase-flow problem, i.e., creeping Stokes flow in the free-flow domain, while using Darcy's law for the porous domain.
While limiting the physical complexity of our model, we focus on the numerical solution of the arising coupled system using either fully monolithic coupled schemes or a partitioned, iterative approach.

We build our contribution on the following observation:
Practitioners, in particular in the modeling community, often rely on sparse direct solvers for the (linearized) subproblems, e.g., \umfpacK{}, \textsc{Pardiso} and \textsc{SuperLU}, see \cite{bollhoefer2019SparseDirectSolvers} for an overview.
This holds when Matlab's Backslash operator or its equivalent in SciPy are used, as they  translate to one of these sparse direct solvers  under the hood.
Often this also applies to users of PDE software frameworks like \dumuX{} \cite{flemisch2011dumux, koch2020dumux}, whose design is in fact aiming to minimize the users' burden of having to deal with every single aspect of the simulation pipeline.

Two issues in this context are often overlooked:
First, sparse direct methods for the linear(ized) system(s) do not scale well in terms of compute time and memory.
Second, the ill-conditioning of a fully assembled monolithic system can lead to severe trustworthiness issues in the solution.
Both issues typically only appear after a model and its corresponding simulation pipeline have been set up, i.e., when test problems are exchanged for real-world scenarios.
Table~\ref{tab: umfpack storage} exemplarily shows the fill-in factors in \umfpacK, when solving the monolithic variant of one of our model problems with the finite volume scheme described in Section~\ref{sec: Model Problem}.
While the matrix density increases mostly linearly due to surface-to-volume arguments, the fill-in for the computed sparse LU decomposition is clearly nonlinear in terms of memory. Thus it translates to compute time for generating the decomposition, and subsequently to solving the linear system using the decomposition.

\begin{table}[htb]
\caption{Memory requirements for storing the sparse matrix and its decomposition measured as matrix entries per degree of freedom.}
\label{tab: umfpack storage}
\centering
\begin{tabular}{lrrrrr}
\toprule
\textbf{DoF} & \textbf{156} & \textbf{1\,056} & \textbf{10\,100} & \textbf{102\,720} & \textbf{1\,001\,000} \\
\textbf{Discretization}  & \textbf{$6\times 6$} & \textbf{$16\times 16$} & \textbf{$50\times 50$} & \textbf{$160\times 160$} & \textbf{$500\times 500$} \\
\midrule
System matrix & 6.7  & 7.3  & 7.6  & 7.7  & 7.7 \\
\umfpack{}   & 19.9 & 38.6 & 78.6 & 145.7& 241.4 \\
\bottomrule
\end{tabular}
\end{table}

In this paper, we demonstrate how carefully devised iterative and thus scalable solvers can alleviate these issues for two different solution strategies for coupled problems:
We consider both a partitioned coupling approach where the subproblems are solved alternately, and a monolithic approach that honors the saddle point structure of the system.
The former is realized with the coupling library preCICE~\cite{preCICE}, while the latter is tailored to standard PDE frameworks like \dumuX.
An important part of our contribution is a thorough comparison of these fundamentally different approaches.

%
%
\section{MODEL PROBLEM} \label{sec: Model Problem}

\begin{figure}[htb]
  \centering
  \captionsetup[subfigure]{singlelinecheck=true}
  \hspace*{-0.4cm}
  \begin{subfigure}[b]{.49\columnwidth}
    \centering
    \footnotesize
    \def\svgwidth{\columnwidth}
    \begingroup\endlinechar=-1
    \resizebox{0.8\columnwidth}{!}{%
    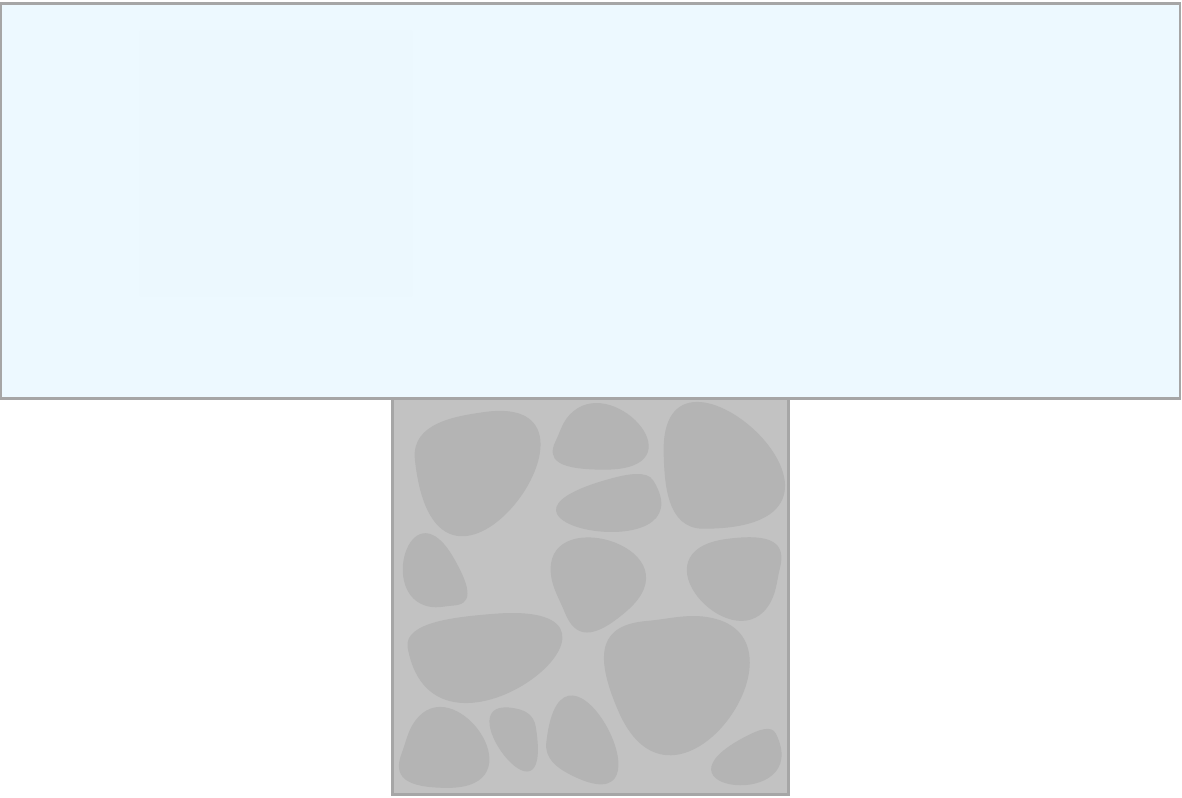%
    }\endgroup
    \caption{Coupled problem illustration.}
    \label{sf: stokes darcy problem setup}
  \end{subfigure}
  ~
  \begin{subfigure}[b]{0.49\columnwidth}
    \centering
    \footnotesize
    \def\svgwidth{\columnwidth}
    \begingroup\endlinechar=-1
    \resizebox{0.8\columnwidth}{!}{%
    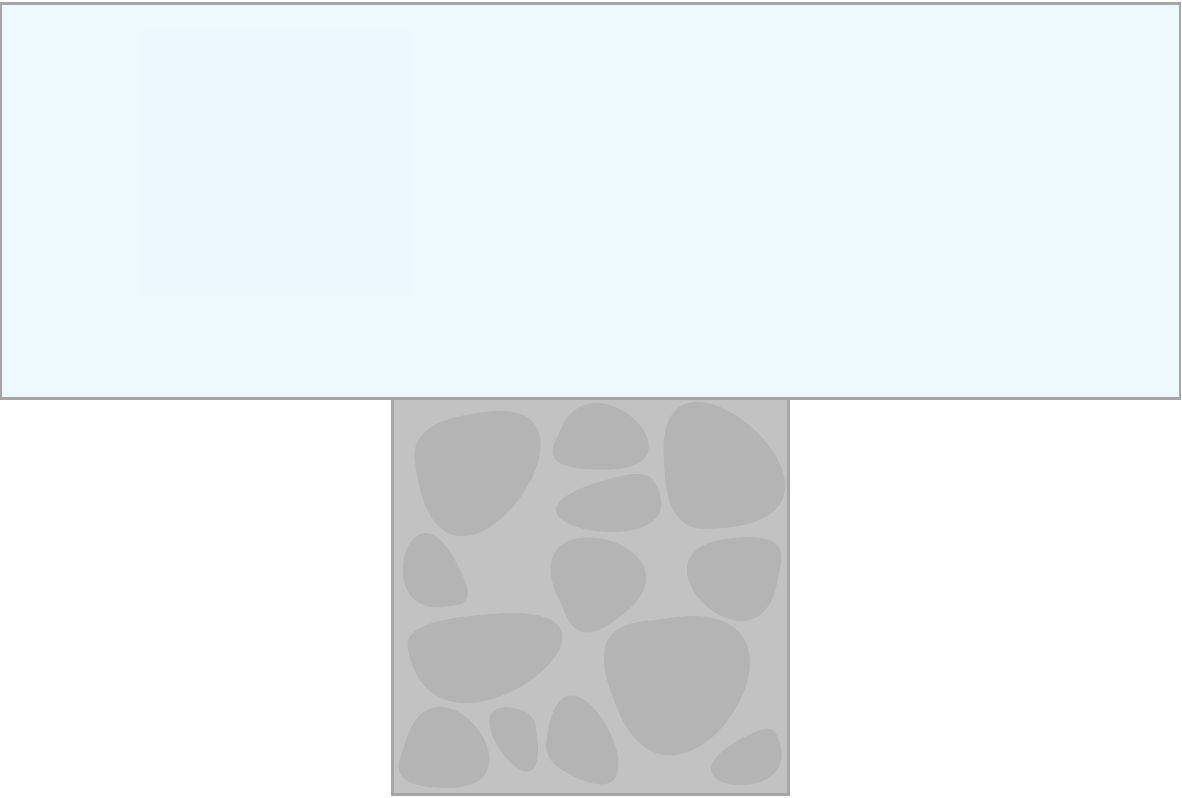%
    }\endgroup
    \caption{Boundary conditions.}
    \label{sf:boundary conditions}
  \end{subfigure}
\vspace{-0.75\baselineskip}
\caption{The instationary Stokes-Darcy problem and its boundary conditions.}
  \label{fig:setup stokes-darcy}
\end{figure}

We consider an instationary, coupled Stokes-Darcy two-domain problem. It comprises a free flow of an incompressible fluid over a porous medium, see
Figure~\ref{sf: stokes darcy problem setup}.

We mark quantities that are associated with the free-flow domain with $\ff$, and use $\prm$ to denote the porous medium.
The domains are denoted by $\Omega\ff$ and $\Omega\prm$, and share the common boundary $\Gamma$.
Boundaries that are not shared are called $\Gamma\ff$ and $\Gamma\prm$.
The normal vectors $n\ff$ and $n\prm$ are orthogonal to $\Gamma$ and point outwards their respective domains.
The time dependent quantities pressure $p$ and velocity $v$ are used to describe the flow in each domain.
We use the transient and incompressible Stokes equations to model the free flow:
\begin{align}
&&&&\frac{\partial v}{\partial t} + \nabla \cdot ( - \nu \left(\nabla v\ff + \nabla v^\mathrm{ff,T} \right)+ \rho^{-1} p\ff I ) &= 0 && \text{in }\Omega\ff&&&&\label{equ: incompress stokes momentum}\\
&&&&\nabla \cdot v\ff & = 0 && \text{in }\Omega\ff&&&&\label{equ: incompress stokes continuity}
\end{align}
Above, $\rho$ and $\nu$ are the fluid density and kinematic viscosity, and $I$ is a suitable identity map.
In the porous medium, Darcy's law
and the continuity equation are used:
\begin{align}
&&&&v\prm & = - K \mu^{-1} \nabla p\prm && \text{in }\Omega\prm&&&&\label{equ: incompress pm darcy}\\
&&&&\nabla \cdot v\prm & = 0 && \text{in }\Omega\prm&&&&\label{equ: incompress pm continuity}
\end{align}
$K$ is the intrinsic permeability of the porous medium and $\mu = \nu \rho$ the  dynamic viscosity of the fluid.
As coupling conditions~\cite{layton2002stokesDarcyExistence}, we use the continuity of the normal stresses~\eqref{equ: incompress coupl cont norm stress}, the Beavers-Joseph-Saffman condition~\cite{saffman1971a} in equation~\eqref{equ: incompress coupl BJS} and the continuity of the normal mass fluxes~\eqref{equ: incompress coupl cont norm mass fluxes}:
\begin{align}
&&&&n \cdot [(p I-\tau)n]\ff & = [p]\prm && \text{on }\Gamma &&&&\label{equ: incompress coupl cont norm stress}\\
&&&&[(v + \sqrt{K}(\alpha_{\mathrm{BJ}}\mu)^{-1} \tau n)\cdot t_\mathrm{ff,pm}]\ff & = 0 && \text{on }\Gamma &&&&\label{equ: incompress coupl BJS}\\
&&&&[v \cdot n]\prm & = -[v \cdot n]\ff && \text{on }\Gamma &&&&\label{equ: incompress coupl cont norm mass fluxes}
\end{align}
We use $n$ for the normal of the respective flow component, $\tau$ for the viscous stresses and $\alpha_{\mathrm{BJ}}$ is the Beavers-Joseph coefficient.
Further, $t_\mathrm{ff,pm}$ is the basis of the tangent plane that describes the interface between $\Omega\ff$ and $\Omega\prm$.
To close this system, boundary conditions for the nonshared domain boundaries are illustrated in Figure~\ref{sf:boundary conditions}.
Note that the pressure $p_\text{in}$ on the left free-flow boundary changes over time.

The system of equations is discretized with a first-order backward Euler scheme in time, and finite volumes in space \cite{koch2020dumux}.
In the Darcy domain, a two-point flux approximation is used for the finite volume approximation of the pressure \cite[Chap.~4]{grueningerChapDiscret}.
In the Stokes domain, a staggered grid is used for the quantities pressure and velocity, and the fluxes are approximated with an upwind scheme \cite{koch2020dumux,schneider2020coupling}.
In summary, the discrete model, to be solved for every time step, has the form
\begin{align}
Ax = b \quad &= \quad \begin{pmatrix} A' & B' \\ C' & D'\end{pmatrix} \begin{pmatrix} x\ff \\ p\prm \end{pmatrix} = \begin{pmatrix} b\ff \\ 0 \end{pmatrix} \label{equ: tp stokesdarcy ff pm}\\
&= \quad \begin{pmatrix} \begin{bmatrix} V & B \\ C & 0 \end{bmatrix} & \begin{bmatrix} B'_1 \\ 0 \end{bmatrix} \\ \begin{bmatrix} C'_1 & 0\end{bmatrix} & D' \end{pmatrix} \begin{pmatrix} \begin{bmatrix} v\ff \\ p\ff \end{bmatrix} \\ p\prm \end{pmatrix} = \begin{pmatrix} \begin{bmatrix} g \\ 0 \end{bmatrix} \\ 0 \end{pmatrix} . \label{equ: stiffness system}
\end{align}
Formulation~\eqref{equ: tp stokesdarcy ff pm} is denoted as \emph{two-domain}~(td) formulation of the problem, because the matrix blocks correspond to the free-flow and porous-medium phase of the problem.
Further, we dub equation~\eqref{equ: stiffness system} the \emph{pressure-velocity}~(pv) formulation, due to the correspondence of the matrix blocks to the variables pressure and velocity.

%
%
\section{PARTITIONED COUPLING APPROACH}

Partitioned coupling approaches are a common strategy to solve coupled problems.
In our setting, this means that the flow fields in the two flow domains are calculated separately, and the coupling between the subdomains is ensured by exchanging information over the sharp interface $\Gamma$.
The benefit of this approach is that existing, optimized solvers for the subdomains can be used. We rely on \dumuX{} for the subdomain solvers, and preCICE for the coupling.

\begin{figure}[htb]
  \centering
  \includegraphics[width=0.63\columnwidth]{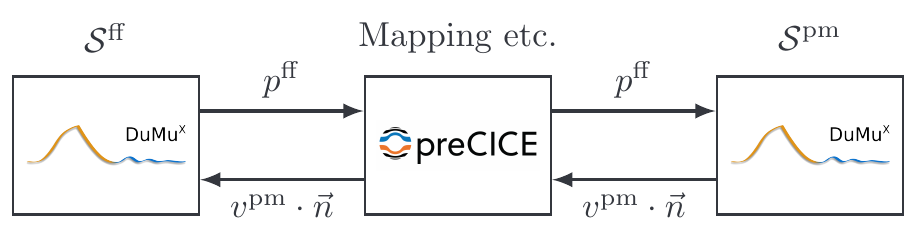}
  \caption{Subdomain coupling scheme implemented with preCICE.}
\label{fig:preCICE coupling hf}
\end{figure}

Looking at the two-domain formulation~\eqref{equ: tp stokesdarcy ff pm},
it is clear that boundary conditions on the common interface $\Gamma$ need to be exchanged in order to get a well-defined solution.
For this, we use a serial implicit coupling technique~\cite{DEGROOTE2013} where the subdomain problems are solved sequentially and the boundary values for the other domain are written after the solution step is performed.
The coupling procedure is depicted in Figure~\ref{fig:preCICE coupling hf} and comprises a Dirichlet-Neumann coupling between the subdomains.
We start the coupling by solving the free-flow problem, to determine Dirichlet pressure values on the interface $\Gamma$.
The porous-medium-flow solver then determines Neumann velocity values on the interface.
Thus, the pressure in the Darcy domain is fixed at the coupling interface $\Gamma$, and the Dirichlet-Neumann coupling leads to a well-defined solution.
In more detail, let $k$ be the coupling iteration index and $v_k^{\mathrm{pm}, \Gamma}$ the normal velocity at the interface $\Gamma$.
The free-flow solver $\solv\ff$ computes a new flow state,
which leads to an updated pressure $p_{k+1}^{\mathrm{ff}, \Gamma}$ on the interface.
With this updated pressure, the porous-medium-flow solver $\solv\prm$ computes a new flow state, which then leads to an updated normal velocity $v_{k+1}^{\mathrm{pm}, \Gamma}$ on the interface.
When we combine the two interface equations
\begin{align}
 \begin{drcases}
 \solv\ff \left(v_k^{\mathrm{pm}, \Gamma}\right) &= p_{k+1}^{\mathrm{ff}, \Gamma} \\
 \solv\prm \left(p_{k+1}^{\mathrm{ff}, \Gamma}\right) &= v_{k+1}^{\mathrm{pm}, \Gamma}
 \end{drcases}
 \quad \Leftrightarrow \quad \solv\prm\left(\solv\ff\left(v_k^{\mathrm{pm}, \Gamma}\right)\right) = v_{k+1}^{\mathrm{pm}, \Gamma},
\end{align}
we can interpret the coupling scheme as an iterative solver for the fixed-point problem
\begin{align}
 \solv\prm\left(\solv\ff\left(v^{\mathrm{pm}, \Gamma}\right)\right) = v^{\mathrm{pm}, \Gamma} \quad \Leftrightarrow \quad R\left(v^{\mathrm{pm}, \Gamma}
 \right) :=  \solv\prm\left(\solv\ff\left(v^{\mathrm{pm}, \Gamma}\right)\right) - v^{\mathrm{pm}, \Gamma} = 0.
 \label{equ:fixed-point problem}
\end{align}
The scheme stops when the interface values converge, i.e., the fixed-point problem is solved to a prescribed accuracy.
We emphasize that when the residual $R$ is sufficiently close to zero, we recover the monolithic solution.

\begin{figure}[htb]
  \centering
  \includegraphics[width=\textwidth]{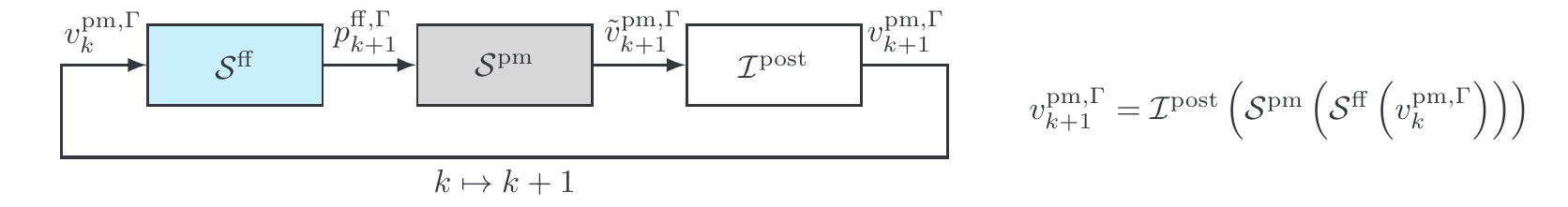}
  \caption{preCICE coupling scheme with enabled post-processing.}
\label{fig:preCICE coupling hf post}
\end{figure}

Solving the fixed-point problem \eqref{equ:fixed-point problem} with a Picard fixed-point iteration is prone to divergence for problems with strong instabilities or oscillations.
In order to improve stability and convergence speed, fixed-point acceleration methods enrich the Picard iteration.
These methods are applied as a post-processing step that we denote as $\postproc$.
Figure~\ref{fig:preCICE coupling hf post} illustrates our accelerated fixed-point iteration.
Now, the flow update from the porous-medium solver is denoted by $\tilde{v}_{k+1}^{\mathrm{pm}, \Gamma}$ (previously: $v_{k+1}^{\mathrm{pm}, \Gamma}$), as it is the solution \emph{before} the improvement by $\postproc$.
The post-processing scheme receives $\tilde{v}_{k+1}^{\mathrm{pm}, \Gamma}$ from the porous-medium solver and computes an improved velocity $v_{k+1}^{\mathrm{pm}, \Gamma}$.
This new velocity depends on the current value and a history of previously calculated  values.
For our experiments, we choose the inverse least-squares interface quasi-Newton method \cite{DEGROOTE2009} for the post-processing.
This method approximates the inverse Jacobian of the residual operator $R$ of the nonlinear coupling equation, based on input-output relations.
Then, it performs Newton-like update steps where a norm minimization is carried out.
For more details on post-processing schemes and their implementation, see \cite{DEGROOTE2013} and \cite{preCICE} respectively.
To determine when the iteration can be stopped, we use the relative convergence measures
\begin{align}
 \left\lVert p_{k+1}^{\mathrm{ff}, \Gamma} - p_k^{\mathrm{ff}, \Gamma} \right\rVert_2 < \varepsilon \left\lVert p_{k+1}^{\mathrm{ff}, \Gamma} \right\rVert_2
 \qquad
 \text{and}
 \qquad
 \left\lVert \tilde{v}_{k+1}^{\mathrm{pm}, \Gamma} - v_k^{\mathrm{pm}, \Gamma} \right\rVert_2 < \varepsilon \left\lVert \tilde{v}_{k+1}^{\mathrm{pm}, \Gamma}
 \right\rVert_2.
\end{align}
Our choice of post-processing method and convergence measure is based on \cite{JAUST2020}.
There, it is shown that for a similar model problem, the coupling approach as outlined above is consistent and converges to the monolithic solution.
This finding is the basis for our convergence study in Section~\ref{RESULTS}.
As solvers $\solv\ff$ and $\solv\prm$, we use problem specific preconditioned iterative subdomain solvers in order to benefit from their smaller memory footprint, which results in a better numerical scaling with respect to the problem size.

preCICE follows a pure library approach, and is called from within \dumux{} by the participating solvers, and not by the general framework.
Through the library approach, the coupling is minimally invasive and provides a black-box functionality that allows us to use optimized solvers for the two domains.

%
%
\section{MONOLITHIC BLOCK PRECONDITIONING}

As an alternative to the partitioned approach, we consider a block-preconditioning strategy to iteratively solve the linear system in its entirety,
allowing to consider both formulations~\eqref{equ: tp stokesdarcy ff pm} and~\eqref{equ: stiffness system} separately.
The saddle point structure of $A'$ implies one of the diagonal matrix blocks to be zero, which prevents the direct application of many preconditioning techniques like simple splitting based schemes~\cite[Chap.~10.2]{saad2003iterative} or incomplete LU (ILU) preconditioning~\cite{meijerink1977iterative,hysom2002ilup}.
To address this, we apply two types of block-preconditioning schemes that use exchangeable preconditioners for the respective matrix blocks.
This allows to select the preconditioners for the matrix blocks based on structural or model-based properties of the block.
The two considered block-preconditioning approaches are a \emph{block-Jacobi} $\pc\bjac$ and a \emph{block-Gauss-Seidel preconditioning scheme} $\pc\bgs$, which themselves are formulated to regard the linear system either as the $2\times 2$ block matrix~\eqref{equ: tp stokesdarcy ff pm} or as the $3\times 3$ block matrix~\eqref{equ: stiffness system}.

To construct the $\pc\bjac$ preconditioners, we consider the block diagonal of the matrices~\eqref{equ: tp stokesdarcy ff pm} and~\eqref{equ: stiffness system}.
Additionally incorporating all block lower triangular parts is the basis for the $\pc\bgs$ preconditioners.
Block-`inverting' those reduced block matrices yields the block preconditioners.
The procedure to acquire the $\pc\bgs$ preconditioner is similar in spirit to \cite{cai_preconditioning_2009}.
There are two points to note:
Firstly, an exact block inversion requires the exact inverses of the diagonal blocks, which is infeasible for preconditioners.
We thus  replace the inverses of the blocks $A'$, $D'$ and $V$ with preconditioners for these blocks that approximate the action of the exact inverses on the quantities of interest.
We denote these block specific preconditioners by $\pc\app$, $\pc\dpp$ and $\pc\vu$.
Secondly, we formally replace the zero block on the diagonal of the reduced matrix~\eqref{equ: stiffness system} with an identity matrix, preventing the `inversion' of the zero block.
In our implementation, this means that no preconditioner is applied to this block.
We thus obtain general $\pc\bjac$ and $\pc\bgs$ preconditioner formulations.

In the implementation, a variety of concrete preconditioners for the respective matrix blocks can be used.
For the two-domain formulation~\eqref{equ: tp stokesdarcy ff pm}, this yields the \emph{two-domain block-Jacobi preconditioner} $\pc\tp\bjac(\pc\app, \pc\dpp)$ and the \emph{two-domain block-Gauss-Seidel preconditioner} $\pc\tp\bgs(\pc\app, \pc\dpp)$.
Both depend on suitable preconditioners $\pc\app$ and $\pc\dpp$ for the blocks $A'$ and $D'$:
\begin{align}
\pc\tp\bjac(\pc\app, \pc\dpp) &\coloneqq  \begin{pmatrix} \pc\app & 0 \\ 0 & \pc\dpp \end{pmatrix}
,&
\pc\tp\bgs(\pc\app, \pc\dpp) &\coloneqq \begin{pmatrix} \pc\app & 0 \\ -\pc\dpp C' \pc\app & \pc\dpp \end{pmatrix}
\end{align}

Likewise, we obtain the \emph{pressure-velocity block-Jacobi preconditioner} $\pc\pv\bjac(\pc\vu, \pc\dpp)$ and the \emph{pressure-velocity block-Gauss-Seidel preconditioner} $\pc\pv\bgs(\pc\vu, \pc\dpp)$.
They depend on preconditioners $\pc\vu$ and $\pc\dpp$ for the blocks $V$ and $D'$:
\begin{align}
\pc\pv\bjac(\pc\vu, \pc\dpp) &\coloneqq \begin{pmatrix}
            \pc\vu        &   0               & 0 \\
           0              &   \text{I}             & 0 \\
             0 & 0 &  \pc\dpp
           \end{pmatrix}
,&
\pc\pv\bgs(\pc\vu, \pc\dpp)
       & \coloneqq \begin{pmatrix}
           \pc\vu                  &   0               & 0 \\
         - \text{I} C \pc\vu             &   \text{I}             & 0 \\
         - \pc\dpp C'_1 \pc\vu & 0 &  \pc\dpp
         \end{pmatrix}
\end{align}

The differences between the four variants directly influence the possible choices of the sub-preconditioners $\pc\app$, $\pc\dpp$ and $\pc\vu$, as well as the computational efficiency of the whole block preconditioner.
We begin with the difference between the two-domain $\pc\tp_{*}$ and pressure-velocity $\pc\pv_{*}$ preconditioner formulations.
The two-domain block preconditioners treat the saddle point structure of $A'$ as one block.
This permits specialized saddle point preconditioning techniques $\pc\app$, like an Uzawa preconditioner \cite{elman_inexact_1994}, which may lead to better results due to their construction for the specific structure.
The pressure-velocity formulation explicitly treats the saddle point structure by introducing an identity preconditioner on the critical diagonal block.
This allows using a wider range of preconditioning techniques for $\pc\vu$ and $\pc\dpp$, like splitting techniques, \ilu{p}
or algebraic multigrid (AMG)~\cite{blatt2012AMGnonsmoothAgg,bastian2012algebraic} preconditioners.
Also, this weakens the prerequisites on the preconditioners and allows using techniques that are not specialized for the problem's structure.

Comparing the $\pc\bjac^{*}$ and $\pc\bgs^{*}$ preconditioners, the sparsity pattern of the block-Jacobi matrices suggests that their application requires fewer computational steps.
The Gauss-Seidel preconditioners use additional coupling entries below the diagonal.
Intuitively, one may expect an improved conditioning with a formulation that makes use of such additional information.
However, this comes at the cost of more computation in the preconditioner application and setup.
Without numerical experiments it is unclear whether the envisioned improvement of the system's condition leads to shorter solution times compared to a preconditioner that is cheaper in its application but less capable to improve the systems condition number.

We highlight that preconditioners are commonly implemented as their application to vectors, i.e., for a vector $x$,  it is of the form $\pc(x)$.
If such an implementation is already given, applying the block preconditioners to a blocked vector is a simple block matrix-vector product.
This allows using existing preconditioner implementations within the blocked preconditioner.
We use the ones available in \dumuX{} for our experiments.

%
%
\section{COMPARISON AND RESULTS} \label{RESULTS}

We now assess the iterative solution of a coupled Stokes-Darcy system with partitioned coupling and block preconditioning, compared to the direct solver  \umfpack \cite{davis1995umfpack}.
We compare the runtime for solving increasingly large systems, and also comment on the memory requirement for all approaches.
In the model problem from Section \ref{sec: Model Problem}, we set $K=10^{-6} \text{m}^2$ and $\alpha_\text{BJ}=1.0$.
The simulation is stopped at $t_\text{end}=50\cdot 10^{5}$s, with a time step size of $\text{d}t=2\cdot 10^{5}$s.
A time dependent pressure difference is applied between the left and right boundary, which changes in the form of a half cosine-wave, with a maximum difference of $10^{-9}$Pa.
The Stokes domain is a $1\times 3$ rectangle over a $1\times 1$ square for the Darcy domain, see also Figure \ref{fig:setup stokes-darcy}.
Each $1\times 1$ square uses the same number of spatial cells.

As solvers, we either use \umfpacK, or preconditioned versions of the iterative solvers PD-GMRES \cite{pdgmresm} or Bi-CGSTAB \cite{bicgstab}.
As preconditioners we use an AMG method \cite{blatt2012AMGnonsmoothAgg,bastian2012algebraic}, Uzawa-iterations \cite{elman_inexact_1994} or an \ilu{0} factorization \cite{hysom2002ilup,meijerink1977iterative}.
PD-GMRES uses $m_\text{init}=m_\text{min}=3$ and $m_\text{step}=5$, other parameters are chosen according to the original publication.
\dunE-\istL's non-smooth aggregation AMG is used as solver or preconditioner, and performs one V-cycle~\cite{blatt2006istl,blatt2012AMGnonsmoothAgg}.
Pre- and post-smoothing is a single Gauss-Seidel iteration each with $\omega=1$.
In our setup, we restrict ourselves to \umfpack{} as coarse grid solver, and limit the hierarchy to 3 levels, which results in comparatively large coarse grid problems for the velocity blocks at scale.
Our standard Uzawa configuration is an inexact Uzawa that executes one Richardson iteration where the optimal relaxation parameter $\omega_{\text{opt}}$ is estimated via power iteration.
In the inexact case, the AMG method as specified above is used as solver, while \umfpack is used for the exact Uzawa iteration ($\text{Uzawa}_e$).
All simulations are run on a single core of an AMD EPYC 7551P CPU with 2.0 GHz.

 \begin{figure}[!tb]
   \centering
   \begin{subfigure}[t]{.25\columnwidth}
   \centering
   \includegraphics[width=\columnwidth]{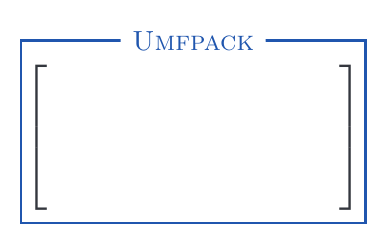}
   \caption{\umfpacK}
   \end{subfigure}
   ~
   \hspace*{-0.03\columnwidth}
   \begin{subfigure}[t]{.24\columnwidth}
   \centering
   \includegraphics[width=\columnwidth]{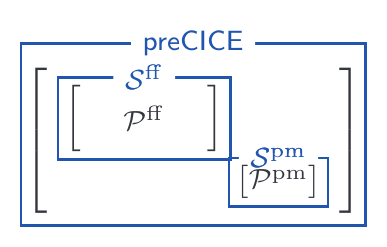}
   \caption{preCICE}
   \end{subfigure}
   ~
   \hspace*{-0.03\columnwidth}
   \begin{subfigure}[t]{.24\columnwidth}
   \centering
   \includegraphics[width=\columnwidth]{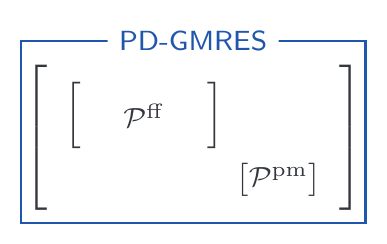}
   \caption{Preconditioner $\pc\tp$}
   \end{subfigure}
   ~
   \hspace*{-0.03\columnwidth}
   \begin{subfigure}[t]{.24\columnwidth}
   \centering
   \includegraphics[width=\columnwidth]{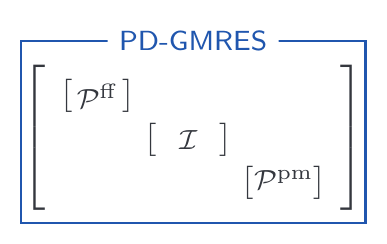}
   \caption{Preconditioner $\pc\pv$}
   \end{subfigure}
   \caption{Evaluation setup conceptual visualization. Solvers are marked in blue around the matrix block they are applied to, preconditioners $\mathcal{P}$ are marked within their block.}
 \label{fig:evaluation concept vis}
 \end{figure}

\begin{table}[tb]
\footnotesize
\centering
\caption{Evaluation setup comparison. Cells marked with \emph{n.a.} are not applicable.}
\label{tab: evaluation concept}
\resizebox{1\textwidth}{!}{%
  \begin{tabular}{lccllllll}
    \toprule
    \textbf{Name} & \textbf{Iterative} & \textbf{Monolithic} & \multicolumn{3}{l}{\textbf{Solver}}  & \multicolumn{3}{l}{\textbf{Preconditioner}} \\
    \cline{4-6} \cline{7-9}
    &  &  & \textbf{System} & \textbf{Stokes} & \textbf{Darcy} & \textbf{System} & \textbf{Stokes} & \textbf{Darcy} \\
    \midrule
    \umfpacK & \no & \yes & \umfpacK{} & \emph{n.a.} & \emph{n.a.} & \emph{n.a.} & \emph{n.a.} & \emph{n.a.}  \\
    \midrule
    preCICE \umfpacK & \yes & \no & preCICE & \umfpacK & \umfpacK & \emph{n.a.} & \emph{n.a.} & \emph{n.a.}\\
    preCICE $\pc(\text{Uzawa}_e, \text{AMG})$  & \yes & \no & preCICE & PD-GMRES & Bi-CGSTAB & \emph{n.a.} & Uzawa-exact & AMG \\
    preCICE $\pc(\text{Uzawa}, \text{AMG})$  & \yes & \no & preCICE & PD-GMRES & Bi-CGSTAB & \emph{n.a.} & Uzawa & AMG \\
    \midrule
    PD-GMRES $\pc\pv\bjac(\text{AMG}, \text{AMG})$ & \yes & \yes & PD-GMRES & \emph{n.a.} & \emph{n.a.} & B-Jacobi & AMG & AMG\\
    PD-GMRES $\pc\pv\bgs(\text{AMG}, \text{AMG})$ & \yes & \yes & PD-GMRES & \emph{n.a.} & \emph{n.a.} & B-Gauss-Seidel  & AMG & AMG\\
    PD-GMRES $\pc\tp\bjac(\text{Uzawa}, \text{AMG})$ & \yes & \yes & PD-GMRES & \emph{n.a.} & \emph{n.a.} & B-Jacobi & Uzawa & AMG \\
    PD-GMRES $\pc\tp\bgs(\text{Uzawa}, \text{AMG})$ & \yes & \yes & PD-GMRES & \emph{n.a.} & \emph{n.a.} & B-Gauss-Seidel & Uzawa & AMG\\
    PD-GMRES $\pc\tp\bjac(\text{Uzawa}, \ilu{0})$ & \yes & \yes & PD-GMRES & \emph{n.a.} & \emph{n.a.} & B-Jacobi & Uzawa & $\ilu{0}$ \\
    PD-GMRES $\pc\tp\bgs(\text{Uzawa}, \ilu{0})$ & \yes & \yes & PD-GMRES & \emph{n.a.} & \emph{n.a.} & B-Gauss-Seidel & Uzawa & $\ilu{0}$\\
    \bottomrule
  \end{tabular}
}
\end{table}
\begin{figure}[thb]
\centering
  \def\svgwidth{0.85\columnwidth}
      \executeiffilenewer{Illustrations/runtimes_paper_hf_fs13.svg}{Illustrations/runtimes_paper_hf_fs13.pdf}%
      {inkscape -z -D --file=Illustrations/runtimes_paper_hf_fs13.svg --export-pdf=Illustrations/runtimes_paper_hf_fs13.pdf --export-latex}%
      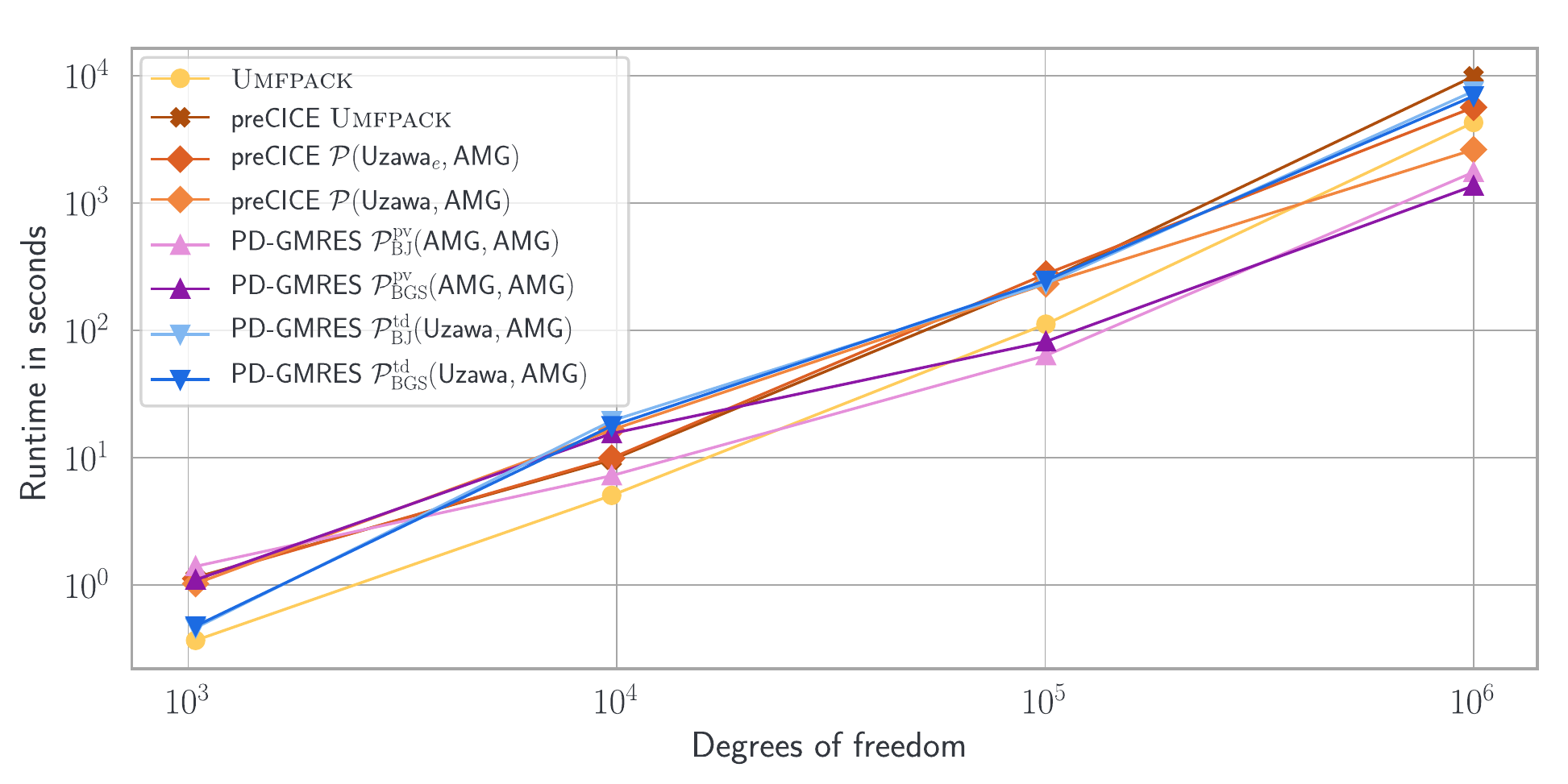%

  \caption{Runtime comparison of the direct solver \umfpack and iterative solvers with partitioned coupling and block preconditioning.}
  \label{fig: Runtime comparison}
\end{figure}

\begin{figure}[thb]
\centering
  \def\svgwidth{0.85\columnwidth}
      \executeiffilenewer{Illustrations/runtimes_paper_bestresults_hf_10004000_fs13.svg}{Illustrations/runtimes_paper_bestresults_hf_10004000_fs13.pdf}%
      {inkscape -z -D --file=Illustrations/runtimes_paper_bestresults_hf_10004000_fs13.svg --export-pdf=Illustrations/runtimes_paper_bestresults_hf_10004000_fs13.pdf --export-latex}%
      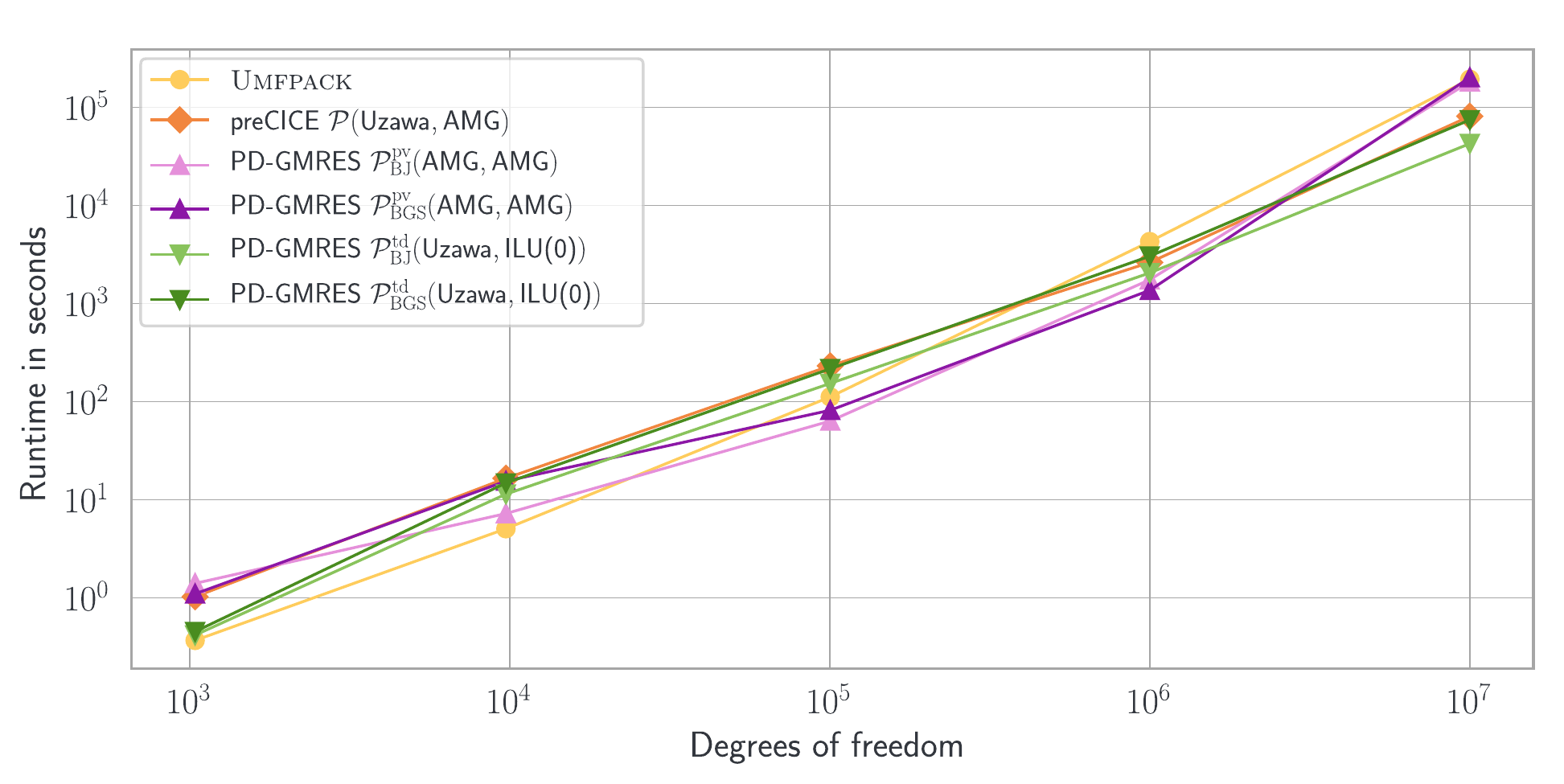%

  \caption{Runtime comparison of the best performing solver configurations.}
  \label{fig: Runtime best}
\end{figure}

To assess the runtime scaling of our different approaches, we increase the number of degrees of freedom.
As baseline, we solve the linear system~\eqref{equ: stiffness system} with the direct solver \umfpacK.
The evaluations tested for the two iterative schemes are listed in Table~\ref{tab: evaluation concept}, and schematically illustrated in Figure~\ref{fig:evaluation concept vis}.
The iterative methods are stopped when the residual's norm is
in the same order of magnitude as the \umfpacK's residual.

Figure~\ref{fig: Runtime comparison} shows the measured runtime scaling behavior.
To allow a comparison between the approaches, we choose the preconditioners for the subsystems to be either AMG or Uzawa. 
We observe that using iterative methods pays off in terms of runtime already for moderate problem sizes with $10^4$ degrees of freedom, benefiting from their better numerical scaling with respect to the problem size $n$.
Partitioned coupling with preconditioned iterative solvers is able to outperform \umfpacK{} for large $n$, while using the partitioned coupling approach with \umfpack for both subsystems is not beneficial and always slower than directly applying \umfpack to the monolithic system.
The performance of our block-preconditioning approach yields the fastest runtimes for large systems, but varies strongly with the preconditioner configuration.
Interestingly, we observe that using the specialized Uzawa preconditioner for Stokes in $\pc\tp_{*}$ leads to increased runtimes compared to the less specialized $\pc\pv_{*}$ block preconditioner with two AMG preconditioners.
In general, the $\pc\bjac^{*}$ configurations lead to slightly improved runtimes compared to the corresponding $\pc\bgs^{*}$ preconditioners.

In Figure \ref{fig: Runtime best} we show that tweaking the preconditioner configurations has the potential to further speed up the runtime, especially for the block-preconditioning approaches. While \umfpack scales roughly as $O(n\cdot\log(n))$, our partitioned and block-preconditioned approaches suggest a linear runtime increase with respect to the problem size $n$. In general, this behavior is also expected for the $\pc\pv_*(\text{AMG}, \text{AMG})$ approaches but due to our setup we see an increase in runtime to a level similar to the \umfpack setting. This is caused by our restriction to use \umfpack as coarse grid solver in the preconditioner, and limiting the multigrid hierarchy to 3 levels: The resulting comparatively large coarse grid problems for the velocity blocks start to dominate the overall runtime. Increasing the multigrid hierarchy and/or switching to more efficient iterative solver for the coarse grid problems is expected to reduce the runtimes to or below the level of our $\pc\tp_*(\text{Uzawa}, \ilu{0})$ approaches.

In terms of memory requirements both considered approaches are very similar when configured to use the same iterative solver(s) and preconditioners.
Then, the memory consumption is dominated by the auxiliary vectors used by the iterative solvers to solve the linear system - in parts or as a whole.
If limited memory is an issue, the partitioned coupling approach has the advantage to solve one subsystem at a time, requiring only the memory for solving the current subsystem.

%
%
\section{CONCLUSION AND FUTURE WORK}

Our experiments clearly indicate that both partitioned coupling and block-preconditioning approaches yield superior performance compared to using a sparse direct solver. This holds already for moderate problems sizes in single-threaded computations, and we expect the benefits to be substantially larger in parallel computations. Our implementation in \dumuX{} is very general, and can in principle also be applied for the nonlinear Navier-Stokes case, or for coupled flows involving more physics. Initial experiments show that our sophisticated coupling/preconditioning techniques are then obligatory, as simple monolithic iterative schemes fail due to the severe ill-conditioning.

%
%
\section{ACKNOWLEDGEMENTS}
This work was financially supported by the German Research Foundation (DFG), within the Collaborative Research Center on Interface-‐Driven Multi-‐Field Processes in Porous Media (SFB 1313, Project Number 327154368).
\printbibliography

\end{document}

%% file: Illustrations/problem_stokesdarcy_paper.pdf_tex
\begingroup%
  \makeatletter%
  \providecommand\color[2][]{%
    \errmessage{(Inkscape) Color is used for the text in Inkscape, but the package 'color.sty' is not loaded}%
    \renewcommand\color[2][]{}%
  }%
  \providecommand\transparent[1]{%
    \errmessage{(Inkscape) Transparency is used (non-zero) for the text in Inkscape, but the package 'transparent.sty' is not loaded}%
    \renewcommand\transparent[1]{}%
  }%
  \providecommand\rotatebox[2]{#2}%
  \newcommand*\fsize{\dimexpr\f@size pt\relax}%
  \newcommand*\lineheight[1]{\fontsize{\fsize}{#1\fsize}\selectfont}%
  \ifx\svgwidth\undefined%
    \setlength{\unitlength}{340.16760606bp}%
    \ifx\svgscale\undefined%
      \relax%
    \else%
      \setlength{\unitlength}{\unitlength * \real{\svgscale}}%
    \fi%
  \else%
    \setlength{\unitlength}{\svgwidth}%
  \fi%
  \global\let\svgwidth\undefined%
  \global\let\svgscale\undefined%
  \makeatother%
  \begin{picture}(1,0.6755631)%
    \lineheight{1}%
    \setlength\tabcolsep{0pt}%
    \put(0,0){\includegraphics[width=\unitlength,page=1]{problem_stokesdarcy_paper.pdf}}%
    \put(0.39565322,0.49726692){\color[rgb]{0.2,0.25882353,0.30196078}\makebox(0,0)[lt]{\lineheight{1.25}\smash{\begin{tabular}[t]{l}free flow $\Omega^\mathrm{ff}$\end{tabular}}}}%
    \put(0.37160958,0.17427277){\color[rgb]{0.2,0.25882353,0.30196078}\makebox(0,0)[lt]{\lineheight{1.25}\smash{\begin{tabular}[t]{l}porous media\end{tabular}}}}%
    \put(0.47955175,0.62422123){\color[rgb]{0.2,0.25882353,0.30196078}\makebox(0,0)[lt]{\lineheight{1.25}\smash{\begin{tabular}[t]{l}$\Gamma^\mathrm{ff}$\end{tabular}}}}%
    \put(0.94176785,0.4953005){\color[rgb]{0.2,0.25882353,0.30196078}\makebox(0,0)[lt]{\lineheight{1.25}\smash{\begin{tabular}[t]{l}$\Gamma^\mathrm{ff}$\end{tabular}}}}%
    \put(0.67719366,0.14798291){\color[rgb]{0.2,0.25882353,0.30196078}\makebox(0,0)[lt]{\lineheight{1.25}\smash{\begin{tabular}[t]{l}$\Gamma^\mathrm{pm}$\end{tabular}}}}%
    \put(0.47073256,0.01009296){\color[rgb]{0.2,0.25882353,0.30196078}\makebox(0,0)[lt]{\lineheight{1.25}\smash{\begin{tabular}[t]{l}$\Gamma^\mathrm{pm}$\end{tabular}}}}%
    \put(0.24671821,0.14798291){\color[rgb]{0.2,0.25882353,0.30196078}\makebox(0,0)[lt]{\lineheight{1.25}\smash{\begin{tabular}[t]{l}$\Gamma^\mathrm{pm}$\end{tabular}}}}%
    \put(0.00860011,0.49673873){\color[rgb]{0.2,0.25882353,0.30196078}\makebox(0,0)[lt]{\lineheight{1.25}\smash{\begin{tabular}[t]{l}$\Gamma^\mathrm{ff}$\end{tabular}}}}%
    \put(0.48638106,0.35424406){\color[rgb]{0,0.31764706,0.61960784}\makebox(0,0)[lt]{\lineheight{1.25}\smash{\begin{tabular}[t]{l}$\Gamma$\end{tabular}}}}%
    \put(0.6514957,0.39121666){\color[rgb]{0.2,0.25882353,0.30196078}\makebox(0,0)[lt]{\lineheight{1.25}\smash{\begin{tabular}[t]{l}$n^\mathrm{pm}$\end{tabular}}}}%
    \put(0.54566554,0.26443465){\color[rgb]{0.2,0.25882353,0.30196078}\makebox(0,0)[lt]{\lineheight{1.25}\smash{\begin{tabular}[t]{l}$n^\mathrm{ff}$\end{tabular}}}}%
    \put(0,0){\includegraphics[width=\unitlength,page=2]{problem_stokesdarcy_paper.pdf}}%
    \put(0.46682693,0.11845269){\color[rgb]{0.2,0.25882353,0.30196078}\makebox(0,0)[lt]{\lineheight{1.25000012}\smash{\begin{tabular}[t]{l}$\Omega^\mathrm{pm}$\end{tabular}}}}%
  \end{picture}%
\endgroup%

%% file: Illustrations/problem_stokesdarcy_paper_bd.pdf_tex
\begingroup%
  \makeatletter%
  \providecommand\color[2][]{%
    \errmessage{(Inkscape) Color is used for the text in Inkscape, but the package 'color.sty' is not loaded}%
    \renewcommand\color[2][]{}%
  }%
  \providecommand\transparent[1]{%
    \errmessage{(Inkscape) Transparency is used (non-zero) for the text in Inkscape, but the package 'transparent.sty' is not loaded}%
    \renewcommand\transparent[1]{}%
  }%
  \providecommand\rotatebox[2]{#2}%
  \newcommand*\fsize{\dimexpr\f@size pt\relax}%
  \newcommand*\lineheight[1]{\fontsize{\fsize}{#1\fsize}\selectfont}%
  \ifx\svgwidth\undefined%
    \setlength{\unitlength}{340.16760606bp}%
    \ifx\svgscale\undefined%
      \relax%
    \else%
      \setlength{\unitlength}{\unitlength * \real{\svgscale}}%
    \fi%
  \else%
    \setlength{\unitlength}{\svgwidth}%
  \fi%
  \global\let\svgwidth\undefined%
  \global\let\svgscale\undefined%
  \makeatother%
  \begin{picture}(1,0.6755631)%
    \lineheight{1}%
    \setlength\tabcolsep{0pt}%
    \put(0,0){\includegraphics[width=\unitlength,page=1]{problem_stokesdarcy_paper_bd.pdf}}%
    \put(0.45750377,0.63304061){\color[rgb]{0.2,0.25882353,0.30196078}\makebox(0,0)[lt]{\lineheight{1.25}\smash{\begin{tabular}[t]{l}$v = 0$\end{tabular}}}}%
    \put(0.82711958,0.4953005){\color[rgb]{0.2,0.25882353,0.30196078}\makebox(0,0)[lt]{\lineheight{1.25}\smash{\begin{tabular}[t]{l}$p = p_\text{out}$\end{tabular}}}}%
    \put(0.46632296,0.01891202){\color[rgb]{0.2,0.25882353,0.30196078}\makebox(0,0)[lt]{\lineheight{1.25}\smash{\begin{tabular}[t]{l}$p=0$\end{tabular}}}}%
    \put(0.3481389,0.14798291){\color[rgb]{0.2,0.25882353,0.30196078}\makebox(0,0)[lt]{\lineheight{1.25}\smash{\begin{tabular}[t]{l}$v = 0$\end{tabular}}}}%
    \put(0.0174193,0.49673873){\color[rgb]{0.2,0.25882353,0.30196078}\makebox(0,0)[lt]{\lineheight{1.25}\smash{\begin{tabular}[t]{l}$p = p_\text{in}(t)$\end{tabular}}}}%
    \put(0.41582754,0.3510964){\color[rgb]{0,0.31764706,0.61960784}\makebox(0,0)[lt]{\lineheight{1.25}\smash{\begin{tabular}[t]{l}Neumann\end{tabular}}}}%
    \put(0,0){\includegraphics[width=\unitlength,page=2]{problem_stokesdarcy_paper_bd.pdf}}%
    \put(0.55060924,0.14798291){\color[rgb]{0.2,0.25882353,0.30196078}\makebox(0,0)[lt]{\lineheight{1.25000012}\smash{\begin{tabular}[t]{l}$v = 0$\end{tabular}}}}%
    \put(0.41790344,0.29378101){\color[rgb]{0,0.31764706,0.61960784}\makebox(0,0)[lt]{\lineheight{1.25000012}\smash{\begin{tabular}[t]{l}Coupling\end{tabular}}}}%
    \put(0.13721584,0.34824856){\color[rgb]{0.2,0.25882353,0.30196078}\makebox(0,0)[lt]{\lineheight{1.25000012}\smash{\begin{tabular}[t]{l}$v = 0$\end{tabular}}}}%
    \put(0.78583739,0.34824856){\color[rgb]{0.2,0.25882353,0.30196078}\makebox(0,0)[lt]{\lineheight{1.25000012}\smash{\begin{tabular}[t]{l}$v = 0$\end{tabular}}}}%
  \end{picture}%
\endgroup%

%% file: Illustrations/runtimes_paper_hf_fs13.pdf_tex
\begingroup%
  \makeatletter%
  \providecommand\color[2][]{%
    \errmessage{(Inkscape) Color is used for the text in Inkscape, but the package 'color.sty' is not loaded}%
    \renewcommand\color[2][]{}%
  }%
  \providecommand\transparent[1]{%
    \errmessage{(Inkscape) Transparency is used (non-zero) for the text in Inkscape, but the package 'transparent.sty' is not loaded}%
    \renewcommand\transparent[1]{}%
  }%
  \providecommand\rotatebox[2]{#2}%
  \newcommand*\fsize{\dimexpr\f@size pt\relax}%
  \newcommand*\lineheight[1]{\fontsize{\fsize}{#1\fsize}\selectfont}%
  \ifx\svgwidth\undefined%
    \setlength{\unitlength}{559.95268196bp}%
    \ifx\svgscale\undefined%
      \relax%
    \else%
      \setlength{\unitlength}{\unitlength * \real{\svgscale}}%
    \fi%
  \else%
    \setlength{\unitlength}{\svgwidth}%
  \fi%
  \global\let\svgwidth\undefined%
  \global\let\svgscale\undefined%
  \makeatother%
  \begin{picture}(1,0.5004014)%
    \lineheight{1}%
    \setlength\tabcolsep{0pt}%
    \put(0,0){\includegraphics[width=\unitlength,page=1]{runtimes_paper_hf_fs13.pdf}}%
  \end{picture}%
\endgroup%

%% file: Illustrations/runtimes_paper_bestresults_hf_10004000_fs13.pdf_tex
\begingroup%
  \makeatletter%
  \providecommand\color[2][]{%
    \errmessage{(Inkscape) Color is used for the text in Inkscape, but the package 'color.sty' is not loaded}%
    \renewcommand\color[2][]{}%
  }%
  \providecommand\transparent[1]{%
    \errmessage{(Inkscape) Transparency is used (non-zero) for the text in Inkscape, but the package 'transparent.sty' is not loaded}%
    \renewcommand\transparent[1]{}%
  }%
  \providecommand\rotatebox[2]{#2}%
  \newcommand*\fsize{\dimexpr\f@size pt\relax}%
  \newcommand*\lineheight[1]{\fontsize{\fsize}{#1\fsize}\selectfont}%
  \ifx\svgwidth\undefined%
    \setlength{\unitlength}{559.95266577bp}%
    \ifx\svgscale\undefined%
      \relax%
    \else%
      \setlength{\unitlength}{\unitlength * \real{\svgscale}}%
    \fi%
  \else%
    \setlength{\unitlength}{\svgwidth}%
  \fi%
  \global\let\svgwidth\undefined%
  \global\let\svgscale\undefined%
  \makeatother%
  \begin{picture}(1,0.50125179)%
    \lineheight{1}%
    \setlength\tabcolsep{0pt}%
    \put(0,0){\includegraphics[width=\unitlength,page=1]{runtimes_paper_bestresults_hf_10004000_fs13.pdf}}%
  \end{picture}%
\endgroup%

%% file: references.bib
@TechReport{bollhoefer2019SparseDirectSolvers,
  author       = {Bollh{\"o}fer, Matthias and Schenk, Olaf and Janal{\'i}k, Radim and Hamm, Steve and Gullapalli, Kiran},
  institution  = {arXiv},
  title        = {State-of-the-Art Sparse Direct Solvers},
  number       = {arXiv:1907.05309},
  howpublished = {\url{https://arxiv.org/abs/1907.05309}},
  year         = {2019},
}

@Article{flemisch2011dumux,
  author    = {B. Flemisch and M. Darcis and K. Erbertseder and B. Faigle and A. Lauser and K. Mosthaf and S. M{\"u}thing and P. Nuske and A. Tatomir and M. Wolff and R. Helmig},
  title     = {{DuMux}: {DUNE} for Multi-\{Phase, Component, Scale, Physics, \textellipsis\} Flow and Transport in Porous Media},
  number    = {9},
  pages     = {1102--1112},
  volume    = {34},
  journal   = {Advances in Water Resources},
  publisher = {Elsevier},
  year      = {2011},
}

@Article{koch2020dumux,
  author    = {Koch, Timo and Gl{\"a}ser, Dennis and Weishaupt, Kilian and Ackermann, Sina and Beck, Martin and Becker, Beatrix and Burbulla, Samuel and Class, Holger and Coltman, Edward and Emmert, Simon and others},
  title     = {{DuMux} 3---An Open-Source Simulator for Solving Flow and Transport Problems in Porous Media with a Focus on Model Coupling},
  journal   = {Computers \& Mathematics with Applications},
  publisher = {Elsevier},
  volume = {81},
  pages = {423--443},
  year = {2021},
}

@article{preCICE,
    title = {{preCICE} -- A fully parallel library for multi-physics surface coupling},
    journal = {Computers and Fluids},
    publisher = {Elsevier},
    volume = {141},
    pages = {250--258},
    year = {2016},
    note = {Advances in Fluid-Structure Interaction},
    issn = {0045-7930},
    doi = {https://doi.org/10.1016/j.compfluid.2016.04.003},
    url = {http://www.sciencedirect.com/science/article/pii/S0045793016300974},
    author = {Bungartz, Hans-Joachim and Lindner, Florian and Gatzhammer, Bernhard and Mehl, Miriam and Scheufele, Klaudius and Shukaev, Alexander and Uekermann, Benjamin}
}

@Article{layton2002stokesDarcyExistence,
  author    = {Layton, William J. and Schieweck, Friedhelm and Yotov, Ivan},
  title     = {Coupling Fluid Flow with Porous Media Flow},
  number    = {6},
  pages     = {2195--2218},
  volume    = {40},
  journal   = {SIAM Journal on Numerical Analysis},
  publisher = {SIAM},
  year      = {2002},
}

@article{saffman1971a,
  title={On the boundary condition at the surface of a porous medium},
  author={Saffman, Philip Geoffrey},
  journal={Studies in Applied Mathematics},
  volume={50},
  number={2},
  pages={93--101},
  year={1971},
  doi={10.1002/sapm197150293},
  publisher={Wiley Online Library}
}

@Book{grueningerChapDiscret,
  author    = {Gr{\"u}ninger, Christoph},
  title     = {Numerical Coupling of Navier--Stokes and Darcy Flow for Soil-Water Evaporation},
  isbn      = {978-3-942036-57-3},
  _pagetotal = {132},
  publisher = {Eigenverlag des Instituts f{\"u}r Wasser- und Umweltsystemmodellierung der Universit{\"a}t Stuttgart},
  url       = {http://elib.uni-stuttgart.de/handle/11682/9674},
  abstract  = {The objective of this work is to develop algorithms and provide a framework for an efficient coupling of free flow and porous-medium flow to simulate porous-medium-soil-water evaporation.
The implementation must particularly be capable of simulating laminar free flows, be fast enough for applied research, and cover simulations in two and three dimensions with complex geometries.},
  file      = {Grüninger - Numerical Coupling of Navier-Stokes and Darcy Flow.pdf:/home/jenny/Documents/Studium/13_Masterarbeit/zotero_MA/storage/NVDQSQ55/Grüninger - Numerical Coupling of Navier-Stokes and Darcy Flow.pdf:application/pdf},
  keywords  = {stokes, navier-stokes, darcy},
  langid    = {english},
  year      = {2017},
}

@Article{schneider2020coupling,
  author    = {Schneider, Martin and Weishaupt, Kilian and Gl{\"a}ser, Dennis and Boon, Wietse M. and Helmig, Rainer},
  title     = {Coupling Staggered-Grid and {MPFA} Finite Volume Methods for Free Flow/Porous-Medium Flow Problems},
  volume    = {401},
  journal   = {Journal of Computational Physics},
  publisher = {Elsevier},
  pages = {109012},
  year      = {2020},
}

@article{DEGROOTE2013,
author = {Degroote, Joris},
year = {2013},
pages = {185--238},
title = {Partitioned Simulation of Fluid-Structure Interaction},
volume = {20},
journal = {Archives of Computational Methods in Engineering},
doi = {10.1007/s11831-013-9085-5}
}

@article{DEGROOTE2009,
title = "Performance of a new partitioned procedure versus a monolithic procedure in fluid-structure interaction",
journal = "Computers \& Structures",
volume = "87",
number = "11",
pages = "793--801",
year = "2009",
note = "Fifth MIT Conference on Computational Fluid and Solid Mechanics",
issn = "0045-7949",
doi = "https://doi.org/10.1016/j.compstruc.2008.11.013",
url = "http://www.sciencedirect.com/science/article/pii/S0045794908002605",
author = "Degroote, Joris and Bathe, Klaus-J{\"u}rgen and Vierendeels, Jan"
}

@InProceedings{JAUST2020,
author = {Jaust, A. and Weishaupt, K. and Mehl, M. and Flemisch, B.},
year = {2020},
pages = {605--613},
title = {Partitioned Coupling Schemes for Free-Flow and Porous-Media Applications with Sharp Interfaces},
booktitle="Finite Volumes for Complex Applications IX - Methods, Theoretical Aspects, Examples",
isbn = {978-3-030-43650-6},
doi = {10.1007/978-3-030-43651-3_57},
publisher = {Springer},
}

@Book{saad2003iterative,
  author    = {Saad, Yousef},
  title     = {Iterative Methods for Sparse Linear Systems},
  edition   = {Second},
  publisher = {SIAM},
  year      = {2003},
}

@Article{meijerink1977iterative,
  author  = {Meijerink, J. A. and van der Vorst, H. A.},
  title   = {An Iterative Solution Method for Linear Systems of Which the Coefficient Matrix Is a Symmetric {M}-matrix},
  pages   = {148--162},
  volume  = {31},
  journal = {Mathematics of Computation},
  year    = {1977},
}

@Article{hysom2002ilup,
  author    = {Hysom, David and Pothen, Alex},
  title     = {Level-Based Incomplete {LU} Factorization: Graph Model and Algorithms},
  journal   = {SIAM Journal on Matrix Analysis and Applications},
  publisher = {SIAM},
  year      = {2002},
}

@Article{cai_preconditioning_2009,
  author       = {Cai, Mingchao and Mu, Mo and Xu, Jinchao},
  journal = {Journal of Computational and Applied Mathematics},
  title        = {Preconditioning Techniques for a Mixed {S}tokes/{D}arcy Model in Porous Media Applications},
  doi          = {10.1016/j.cam.2009.07.029},
  issn         = {0377-0427},
  number       = {2},
  pages        = {346--355},
  url          = {http://www.sciencedirect.com/science/article/pii/S0377042709004269},
  urldate      = {2019-11-27},
  volume       = {233},
  abstract     = {We study numerical methods for a mixed Stokes/Darcy model in porous media applications. The global model is composed of two different submodels in a fluid region and a porous media region, coupled through a set of interface conditions. The weak formulation of the coupled model is of a saddle point type. The mixed finite element discretization applied to the saddle point problem leads to a coupled, indefinite, and nonsymmetric linear system of algebraic equations. We apply the preconditioned {GMRES} method to solve the discrete system and are particularly interested in efficient and effective decoupled preconditioning techniques. Several decoupled preconditioners are proposed. Theoretical analysis and numerical experiments show the effectiveness and efficiency of the preconditioners. Effects of physical parameters on the convergence performance are also investigated.},
  file         = {Preconditioning techniques for a mixed Stokes-Darcy model in porous media applications.pdf:/home/jenny/Documents/Studium/13_Masterarbeit/zotero_MA/storage/8PKIVDDF/Cai et al. - 2009 - Preconditioning techniques for a mixed StokesDarc.pdf:application/pdf},
  keywords     = {preconditioner, stokes, darcy, decoupled, {GMRes}},
  langid       = {english},
  shortjournal = {Journal of Computational and Applied Mathematics},
  year         = {2009},
  _month        = {11},
  day          = {15},
}

@Article{elman_inexact_1994,
  author       = {Elman, Howard C. and Golub, Gene H.},
  journal = {{SIAM} Journal on Numerical Analysis},
  title        = {Inexact and Preconditioned Uzawa Algorithms for Saddle Point Problems},
  doi          = {10.1137/0731085},
  issn         = {0036-1429, 1095-7170},
  number       = {6},
  pages        = {1645--1661},
  url          = {http://epubs.siam.org/doi/10.1137/0731085},
  urldate      = {2019-11-25},
  volume       = {31},
  langid       = {english},
  shortjournal = {{SIAM} J. Numer. Anal.},
  year         = {1994},
}

@techreport{blatt2012AMGnonsmoothAgg,
  author  = {Blatt, Markus and Ippisch, Olaf and Bastian, Peter},
  title   = {A Massively Parallel Algebraic Multigrid Preconditioner Based on Aggregation for Elliptic Problems with Heterogeneous Coefficients},
  institution = {arXiv},
  number = {arXiv:1209.0960},
  year    = {2012},
}

@article{bastian2012algebraic,
  title={Algebraic multigrid for discontinuous Galerkin discretizations of heterogeneous elliptic problems},
  author={Bastian, Peter and Blatt, Markus and Scheichl, Robert},
  journal={Numerical Linear Algebra with Applications},
  volume={19},
  number={2},
  pages={367--388},
  year={2012},
  publisher={Wiley Online Library}
}

@Article{davis1995umfpack,
  author  = {Davis, Timothy A.},
  title   = {{UMFPACK} User Guide},
  journal = {Website: http://www.suitesparse.com},
  year    = {2018},
}

@Article{pdgmresm,
  author  = {Rolando Cuevas N{\'u}{\~n}ez and Christian E. Schaerer and Amit Bhaya},
  title   = {A Proportional-Derivative Control Strategy for Restarting the {GMRES(m)} Algorithm},
  doi     = {j.cam.2018.01.009},
  pages   = {209--224},
  volume  = {337},
  journal = {Journal of Computational and Applied Mathematics},
  year    = {2018},
}

@article{bicgstab,
  title={{B}i-{CGSTAB}: A fast and smoothly converging variant of {B}i-{CG} for the solution of nonsymmetric linear systems},
  author={van der Vorst, H. A.},
  journal={SIAM Journal on scientific and Statistical Computing},
  volume={13},
  number={2},
  pages={631--644},
  year={1992},
  publisher={SIAM}
}

@InProceedings{blatt2006istl,
  author       = {Blatt, Markus and Bastian, Peter},
  booktitle    = {International Workshop on Applied Parallel Computing},
  title        = {The Iterative Solver Template Library},
  organization = {Springer},
  pages        = {666--675},
  year         = {2006},
}
